\documentclass{article}
\usepackage{amssymb}
\usepackage{amsmath}

\usepackage{cite}
\usepackage{amssymb}
\usepackage{amsfonts}
\usepackage{amscd}
\usepackage[matrix]{xypic}
\newtheorem{theorem}{Theorem}

\newtheorem{definition}[theorem]{Definition}

\newtheorem{proposition}[theorem]{Proposition}

\newcommand{\A}{\mathcal{A}}
\newcommand{\B}{\mathcal{B}}
\newcommand{\p}{\mathcal{P}}

\newcommand{\K}{\Bbb{K}}
\newcommand{\AB}{\mathcal{A}\otimes \mathcal{B}}

\begin{document}

\title{\bf On the algebras obtained by tensor product. }
\author{Elisabeth REMM \thanks{%
corresponding author: e-mail: E.Remm@uha.fr} - Michel
GOZE \thanks{M.Goze@uha.fr.}\\
\\
{\small Universit\'{e} de Haute Alsace, F.S.T.}\\
{\small 4, rue des Fr\`{e}res Lumi\`{e}re - 68093 MULHOUSE - France}}
\date{}
\maketitle


\bigskip

\noindent {\bf Keywords.} 
 Tensor product of algebras on quadratic algebras.

\begin{abstract}
Let $\cal{P}$ be a quadratic operad. We define an associated operad $\tilde {\cal{P}}$ such that
for any $\cal{P}$-algebra $\A$ and $\tilde {\cal{P}}$-algebra $\B,$ the algebra  $\A\otimes \B$ is always
 a $\cal{P}$-algebra for the classical tensor product. 
\end{abstract}

\section{Introduction}

Let $\mathcal{P}$ be a quadratic operad with one generating operation (i.e. the algebras on this operad have only one 
operation) and $\mathcal{P}^!$ its dual operad. It satisfies
$\mathcal{P}^!=hom(\mathcal{P},\mathcal{L}ie)$ where $\mathcal{L}ie$ is the quadratic operad corresponding to Lie algebras.
For any $\mathcal{P}$-algebra $\A$ and $\mathcal{P}^!$-algebra $\B$, the vector space $\A \otimes \B$ 
is naturally provided with a Lie algebra product 
\begin{eqnarray}
\label{pdt-tens}
\mu(a_1 \otimes b_1,a_2 \otimes b_2)=\mu_\A(a_1,a_2) \otimes  \mu_\B(b_1,b_2) -\mu_\A(a_2,a_1) \otimes  \mu_\B(b_2,b_1).
\end{eqnarray}
where $\mu_\A$ (resp. $\mu_\B$) is the multiplication of $\A$ (resp. $\B$).
We deduce that the "natural" tensor product   
$ \mu_{\A \otimes \B} = \mu_{\A} \otimes \mu_{\B}$ provides $\A \otimes \B$ with a Lie-admissible algebra structure. 

In \cite{G.R} we have defined special classes of Lie-admissible algebras with relations of definition determined by
an action of the subgroups $G_i$ of the 3-degree symmetric group $\Sigma_3$. 
We obtain quadratic operads, denoted by $G_i-\mathcal{A}ss$ and in this family  we find operads of
Lie-admissible, associative, Vinberg and pre-Lie algebras. For these operads we have proved that for every 
$\mathcal{P}$-algebra $\A$ and $\mathcal{P}^!$-algebra $\B$ the tensor product $\A \otimes \B$ is a $\mathcal{P}$-algebra.
This is not true for general nonassociative algebras and, in this sense, the $G_i$-associative algebras are the
most regular kind of nonassociative algebras. For example if $\mathcal{P}$ is the operad of Leibniz algebras or 
of the nonassociative algebras associated to Poisson algebras \cite{G.R3}, then the tensor product of a 
$\mathcal{P}$-algebra  and $\mathcal{P}^!$-algebra is not a 
$\mathcal{P}$-algebra. 

So we introduce a quadratic operad, denoted by $\tilde{ \mathcal{P}}$, such that 
the tensor product of a $\mathcal{P}$-algebra with a $\tilde{ \mathcal{P}}$-algebra is a 
$\mathcal{P}$-algebra. In case of $\mathcal{P}=\mathcal{L}ie$ or $G_i$-$ \mathcal{A}ss$ then 
$\tilde{ \mathcal{P}}=\mathcal{P}^!$ this explain the above remarks.         

\section{Nonassociative algebras and operads}

 We assume in this work that $\mathcal{P}$ is a quadratic operad 
with one generating operation.  Then it is defined from
 the free operad $\Gamma (E)$ generated by a $\Sigma_2$-module 
$E$ placed in arity 2, and an ideal $(R)$ generated by a $\Sigma_3$-invariant subspace $R$ of $\Gamma (E)(3).$
Our hypothesis implies that the $\Sigma_2$-module $E$ is generated by one element (i.e. algebras over this operad are
 algebras with one operation).
Recall that, if we consider an operation with no symmetry, the $\mathbb{K}$-vector space $\Gamma (E)(2)$ is $2$-dimensional
with basis $\{ x_1 \cdot x_2, x_2 \cdot x_1   \}$ and $\Gamma (E)(3)$ is the $12$-dimensional $\mathbb{K}$-vector space generated by 
$\{ x_i \cdot (x_j \cdot x_k), (x_i \cdot x_j) \cdot x_k   \}$ with $i \neq j\neq k \neq i, \, i,j,k\in \{ 1,2,3 \}$.
We have a natural action of $\Sigma_3$  on $\Gamma (E)(3)$ given by:
$$
\begin{array}{llll}
\Sigma_3 &\times & \Gamma (E)(3) & \rightarrow \Gamma (E)(3) \\
( \sigma & , & X ) & \mapsto \sigma (X)
\end{array}  
$$  
where 
$$
\begin{array}{l}
\sigma (x_i \cdot (x_j \cdot x_k)) =x_{\sigma ^{-1}(i)} \cdot (x_{\sigma ^{-1}(j)} \cdot x_{\sigma ^{-1}(k) })\\
\sigma ((x_i \cdot x_j) \cdot x_k )=(x_{\sigma ^{-1}(i)} \cdot x_{\sigma ^{-1}(j)}) \cdot x_{\sigma ^{-1}(k) }).
\end{array}
$$
We denote by $\mathcal{O}(X)$ the orbit of $X$ associated with this action and by $\mathbb{K}(\mathcal{O}(X))$ 
the $\Sigma_3$-invariant subspace of $\Gamma (E)(3)$ generated by $\mathcal{O}(X)$. More generally, if
$X_1,\cdots,X_k$ are vectors in $\Gamma (E)(3)$, we denote
by $\K(\mathcal{O}(X_1,..,X_k))$ the  $\Sigma_3$-invariant subspace of $\Gamma (E)(3)$ generated by 
$\mathcal{O}(X_1)\cup \cdots \cup \mathcal{O}(X_k)$. 

\begin{definition} We say that a $\Sigma_3$-invariant subspace $F$ of $\Gamma (E)(3)$ is of rank $k$ if there exists 
$X_1, ... ,X_k\in \Gamma (E)(3)$ linearly independent such that $F=\mathbb{K}(\mathcal{O}(X_1,...,X_k))$ and 
 for every $ p<k$ and $  Y_1,...Y_p\in F$ we have $\mathbb{K}(\mathcal{O}(Y_1,...,Y_p)) \neq F.$ 
\end{definition}

If $\mathcal{P}$ is a quadratic operad with one generating operation, its module  of relations $R$ is an invariant
subspace of $\Gamma (E)(3)$. We will say that $\mathcal{P}$ is of rank $k$ if and only if $R$ is of rank $k$.

\medskip

\noindent From the action of $\Sigma _3$ on $\Gamma (E)(3)$ we define linear maps on this module as follows. 
Let $\K[\Sigma _3]$ be  the group algebra of
$\Sigma _3$ that is the vector space of all finite linear combinations of elements of $\Sigma _3$ 
with coefficients in $\K$ hence of all elements of the form
$v=a_1\sigma _1+a_2\sigma _2+...+a_6\sigma _6$. If $v=\Sigma a_i\sigma _i \in \K[\Sigma _3]$, let $\Psi _v$ be
given by
$$\Psi _v(X)=\Sigma a_i\sigma _i(X).$$
Then an invariant subspace $F$ of $\Gamma (E)(3)$ is stable for every $\Psi _v$.

\medskip

\noindent Let $(\mathcal{A},\mu)$ be a $\p$-algebra. This means that $(\mathcal{A},\mu)$
is a nonassociative algebra (by nonassociative algebra we mean algebras with 
non necessarily associative multiplication).  We
consider the maps $A^L(\mu)=\mu \circ (\mu \otimes Id)$ and $A^R(\mu)=\mu \circ (Id \otimes \mu).$ 
Then the associator of $\mu$ is written $A(\mu)=A^L(\mu)-A^R(\mu).$
For each vector $v \in \K[\Sigma _3]$ we define
a linear map on $\A^{\otimes 3}$ denoted by $\Phi _v^{\A}$ and given by 
$$
\begin{array}{rcl}
\Phi_v ^\A:\A^{\otimes 3 } & \rightarrow & \A^{\otimes 3 } \\
(x_1 \otimes x_2 \otimes x_3) & \mapsto & \sum a_i (x_{\sigma^{-1}(1)} \otimes x_{\sigma^{-1}(2)}\otimes x_{\sigma^{-1}(3)}) 
\end{array}
$$
The multiplication $\mu$ satisfies relations of the following type
\begin{equation}
\label{nonass}
A^L(\mu) \circ \Phi_v^{\A}- A^R(\mu) \circ \Phi_{v'}^{\A}=0,
\end{equation}
where $v=\Sigma a_i\sigma _i,v'=\Sigma a'_i\sigma _i \in \mathbb{K}[\Sigma_3].$

Such a relation defines the module $R$ of relations of $\mathcal{P}.$ In fact $R$ is the $\Sigma _3$-submodule of
$\Gamma (E)(3)$ is generated by the vectors
$$\sigma _j(\Psi _v(x_1\cdot (x_2\cdot x_3))-\Psi _{v'}((x_1\cdot x_2)\cdot x_3))$$
for every $\sigma _j \in \Sigma _3$.

\medskip

\medskip

Let us note that if $\mathcal{P}$ is of
rank $1$, a $\mathcal{P}$-algebra is given by one multiplication satisfying only one relation of type (\ref{nonass}).

\medskip

\begin{proposition}
Let $\mathcal{P}$ be a quadratic operad with one generating operation such that the $\Sigma _3$-submodule $R$ 
of relations is generated by
vectors of the following type :
$$\Sigma_{i=1}^6 a_i^l\sigma _i((x_1\cdot x_2)\cdot x_3-x_1\cdot (x_2\cdot x_3))$$
for $l=1,...,k$ and $\sigma _i \in \Sigma _3$. Then $\mathcal{P}$ is of rank $1$.
\end{proposition}
{\it Proof.} In fact, the $\Sigma _3$-invariant subspace of $\Gamma (E)(3)$ generated by 
$$((x_1\cdot x_2)\cdot x_3-x_1\cdot (x_2\cdot x_3))$$
is isomorphic to $\K[\Sigma _3]$. This isomorphism is given by:
$$\Sigma a_i\sigma _i((x_1\cdot x_2)\cdot x_3-x_1\cdot (x_2\cdot x_3))\longrightarrow \Sigma a_i\sigma _i.$$
We  have seen in \cite{G.R2} that for every 
$\Sigma _3$-invariant subspace $F$ of $\K[\Sigma _3]$, there is a vector $v \in \K[\Sigma _3]$ such that 
$F=F_v=\K(\mathcal{O}(v))$ where $\mathcal{O}(v)$ is the orbit of $v$ corresponding
to the natural action of $\Sigma _3$ on $\K[\Sigma _3]$. We deduce that the rank is $1$.

\medskip

In the following examples we will recall the definition of the operads $G_i$-$\mathcal{A}ss$ and define some quadratic operads 
of rank one associated to some classes of nonassociative algebras. 

\medskip

\noindent {\bf Examples.} 

\medskip

\noindent {\bf 1. The operads  $G_i\!-\!\mathcal{A}ss$}

\medskip

\noindent 
For $i,j,k \in \{1,2,3 \}$ and $i \neq j \neq k \neq i$ we denote by $\tau_{ij}$ the transposition
 $(i,j)$ and $c_1,c_2$ the cycles $(1,2,3)$ and $(1,3,2)$. 
The subgroups of $\Sigma_3$ are $G_1=\{Id\}, G_2=<\tau_{12}>,G_3=<\tau_{23}>,
G_4=<\tau_{13}>$,
$G_5=<c_1>$, $G_6=\Sigma_3$ where $<g>$ denotes the subgroup generated by $g.$ 
Each one of these subgroups $G_i$ defines 
an invariant submodule $R_i$ of $\Gamma (E)(3)$ of rank 1. In fact consider the vector 
$X=x_1 \cdot(x_2 \cdot x_3)-(x_1 \cdot x_2) \cdot x_3$ of $\Gamma (E)(3)$ and if $G_i$ is a subgroup 
of $\Sigma_3$, we define the vector $X_i$ of $\Gamma (E)(3)$ by
$$X_i=\sum_{\sigma \in G_i}(-1)^{\epsilon(\sigma)}\sigma(X)$$
where $\epsilon(\sigma)$ is the sign of the permutation $\sigma$. 
Let $R_i$ be the subspace $R_i=\mathbb{K}\mathcal({O}(V_i)).$ Then

$$\begin{array}{ll}
R_1= & Vect_{\mathbb{K}}\left\{  (x_i \cdot x_j) \cdot x_k-x_i \cdot (x_j \cdot x_k) \right\} ,\\
R_2= & Vect_{\mathbb{K}}\left\{  (x_i \cdot x_j) \cdot x_k-x_i \cdot (x_j \cdot x_k)-(x_j \cdot x_i) 
\cdot x_k+x_j \cdot (x_i \cdot x_k)  \right\},\\
R_3= & Vect_{\mathbb{K}}\left\{  (x_i \cdot x_j) \cdot x_k-x_i \cdot (x_j \cdot x_k)-(x_i \cdot x_k) 
\cdot x_j+x_i \cdot (x_k \cdot x_j)  \right\},\\ 
R_4= & Vect_{\mathbb{K}}\left\{  (x_i \cdot x_j) \cdot x_k-x_i \cdot (x_j \cdot x_k)-(x_k \cdot x_j) 
\cdot x_i+x_k \cdot (x_j \cdot x_i) \right\}\\
R_5  = & Vect_{\mathbb{K}} \{  (x_i \cdot x_j) \cdot x_k-x_i \cdot (x_j \cdot x_k) +
(x_j \cdot x_k) \cdot x_i-x_j \cdot (x_k \cdot x_i)\\
 &+(x_k \cdot x_i) \cdot x_j-x_k \cdot (x_i \cdot x_j) \\
R_6  = & Vect_{\mathbb{K}} \{  (x_1 \cdot x_2) \cdot x_3-x_1 \cdot (x_2 \cdot x_3) + (x_2 \cdot x_3) \cdot x_1-x_2 \cdot (x_3 \cdot x_1) \\
& +(x_3 \cdot x_1) \cdot x_2-x_3 \cdot (x_1 \cdot x_2) 
  -(x_2 \cdot x_1) \cdot x_3+x_2 \cdot (x_1 \cdot x_3) -(x_3 \cdot x_2) \cdot x_1 \\
&+x_3 \cdot (x_2 \cdot x_1) -(x_1 \cdot x_3) \cdot x_2+x_1 \cdot (x_3 \cdot x_2) \}
\end{array}$$

\begin{definition}
The quadratic operad $G_i$-$\mathcal{A}ss$ is the quadratic operad $\Gamma (E) / (R)$ where $(R)(3)=R_i.$
\end{definition}
Some of these operads are wellknown:

- $G_1$-$\mathcal{A}ss=\mathcal{A}ss$,

- $G_2$-$\mathcal{A}ss=\mathcal{V}inb$,

- $G_3$-$\mathcal{A}ss=\mathcal{P}re-Lie.$

\noindent Let us note that the $G_6$-$\mathcal{A}ss$-algebras are the Lie-admissible algebras that is if  $\mu$
 is the product of such an algebra then
$$[x,y]=\mu(x,y)-\mu(y,x)$$
is a product of Lie algebra. As $G_i$ is a subgroup of $G_6=\Sigma_3$ any $G_i$-$\mathcal{A}ss$-algebra is Lie-admissible. 
For a general study of the operad $G_6$-$\mathcal{A}ss=\mathcal{L}ie\mathcal{A}dm$ and $G_2$-$\mathcal{A}ss=\mathcal{V}inb$ see 
\cite{G.R}.

Let $(G_i$-$\mathcal{A}ss)^!$ be the quadratic dual operad of $G_i$-$\mathcal{A}ss$. 
We denote by $(R_i)^!$ the submodule of $\Gamma (E)(3)$ defining $(G_i$-$\mathcal{A}ss)^!.$
It is the orthogonal of $R_i$ with respect to the inner product on $\Gamma (E)(3)$ given by
\begin{eqnarray*}
&\left<i\cdot (j\cdot k),i\cdot (j\cdot k)\right>=\epsilon\left(
\begin{tabular}{lll}
$1$ & $2$ & $3$ \\
$i$ & $j$ & $k$%
\end{tabular}
\right), \
&\left<(i \cdot j)\cdot k,(i\cdot j)\cdot k\right>=-\epsilon\left(
\begin{tabular}{lll}
$1$ & $2$ & $3$ \\
$i$ & $j$ & $k$%
\end{tabular}
\right)
\end{eqnarray*}
where $\epsilon(\sigma)$ is the sign of the permutation $ \sigma$.
Then we have
$$\begin{array}{ll}
(R_1)^! = & R_1 \\
(R_2)^! = & Vect_{\mathbb{K}}\left\{  (x_i \cdot x_j) \cdot x_k-x_i \cdot (x_j \cdot x_k);(x_i \cdot x_j) 
\cdot x_k- (x_j \cdot x_i)\cdot  x_k \right\}\\
(R_3)^! = & Vect_{\mathbb{K}}\left\{  (x_i \cdot x_j) \cdot x_k-x_i \cdot (x_j \cdot x_k);
(x_i \cdot x_j )\cdot x_k)-(x_i \cdot x_k) \cdot x_j  \right\} \\
(R_4)^! = & Vect_{\mathbb{K}}\left\{  (x_i \cdot x_j) \cdot x_k-x_i \cdot (x_j \cdot x_k);
(x_i \cdot x_j) \cdot x_k-(x_k \cdot x_j) \cdot x_i \right\}\\
(R_5)^! = & Vect_{\mathbb{K}}\left\{  (x_i \cdot x_j) \cdot x_k-x_i \cdot (x_j \cdot x_k);
(x_i \cdot x_j) \cdot x_k-(x_j \cdot x_k) \cdot x_i \right\}\\
(R_6)^! = & Vect_{\mathbb{K}} \{  (x_i \cdot x_j) \cdot x_k-x_i \cdot (x_j \cdot x_k) ;
(x_i \cdot x_j) \cdot x_k -(x_j \cdot x_i) \cdot x_k\\                                          
&(x_i \cdot x_j) \cdot x_k-(x_i \cdot x_k)\cdot x_j \}
\end{array}$$
\begin{proposition} 
For $i=1$, the operad $(G_1$-$\mathcal{A}ss)^!=\mathcal{A}ss$ is if rank $1$. For $2\leq i\leq 6$, the operads
$(G_i$-$\mathcal{A}ss)^!$ are of rank $2$.
\end{proposition}
{\it Proof.} 
The case $i=1$ is trivial (it is also a consequence of Proposition 2).

\noindent For $i=2,3,4$ and $5$, the rank of $(R_i)^!$ is 2. We denote by $v_i^j, j=1,2$ 
the generators  of $(R_i)^!.$ Then if $\B$ is a $(G_i$-$\mathcal{A}ss)^!$-algebra the multiplication $\mu_\B$ satisfies
$$\begin{array}{ll}
1) \ \mu_\B(\mu_\B \otimes Id) -\mu_\B(Id \otimes \mu_\B)=0 \\
2) \ \mu_\B(\mu_\B \otimes Id) -\mu_\B(Id \otimes \mu_\B \circ \Phi _{\sigma_i})=0
\end{array}$$
with 
$$\left\{ 
\begin{array}{l}
\sigma_2=\tau_{12} \\
\sigma_3=\tau_{23} \\
\sigma_4=\tau_{13} \\
\sigma_5=c_1 \ \mbox{\rm{ or }} \ c_2 \\
\end{array}
\right.$$
For $i=6$, the space $(R_6)^!$  is generated by  the vectors 
$$  (x_i \cdot x_j) \cdot x_k-x_i \cdot (x_j \cdot x_k), \ 
(x_i \cdot x_j) \cdot x_k -(x_j \cdot x_i) \cdot x_k, \                                           
(x_i \cdot x_j) \cdot x_k-(x_i \cdot x_k)\cdot x_j.$$ But we can write 
$$(x_i \cdot x_j) \cdot x_k -(x_j \cdot x_i) \cdot x_k= (Id-\tau _{ij})((x_i \cdot x_j)\cdot x_k)$$
and
 $$(x_i \cdot x_j) \cdot x_k -(x_i \cdot x_k) \cdot x_j= (Id-\tau _{jk})((x_i \cdot x_j)\cdot x_k).$$
The $\Sigma _3$-invariant subspace of $\K[\Sigma _3]$ generated by the vectors $Id-\tau _{12}$ and $Id-\tau _{23}$
is of dimension $5$, and from the classification \cite{G.R2}, this space corresponds to $F_v=\K(\mathcal{O})(v)$
with 
$$v=2Id- \tau _{12}-\tau _{13}-\tau _{23}+c_1$$
and we deduce that this operad is of rank $2$.

\medskip

\noindent {\bf 2. The  $3$-power associative algebras}

\medskip

\noindent We have seen that every $(G_i$-$\mathcal{A}ss)$-algebra is Lie-admissible. Moreover the operad 
$\mathcal{L}ie\mathcal{A}dm$
is quadratic, of rank 1. The submodule $R_6$ is of dimension $1$ and corresponds to the one-dimensional
$\Sigma_3$-invariant subspace $\mathbb{K}(\mathcal{O}(V))=F_V$ of $\mathbb{K}[ \Sigma_3],$ where $V$ is given  by
$$V=\sum_{\sigma \in \Sigma_3}(-1)^{\epsilon(\sigma)}\sigma.$$ If we consider the natural action of $\Sigma_3$
on $\mathbb{K} [\Sigma_3],$ then there exists only two irreducible invariant one dimensional subspaces 
of $\mathbb{K}[ \Sigma_3]$ that is $F_V$ and $F_W$ where 
$$W= \sum_{\sigma \in \Sigma_3 } \sigma.$$
 If the set of
Lie-admissible is associated to $F_V$, the set of algebras corresponding to $F_W$ is the set of 
$3$-power associative algebras (see  \cite{G.R2}) that is which satisfies $x^2\cdot x=x\cdot x^2$ for every $x$. 
As in the Lie-admissible case we can define 
classes of $3$-power associative algebras corresponding to the action of the subgroups of $\Sigma_3$.
We are conduced to consider the following submodules of $\Gamma (E)(3)$:
$$\begin{array}{ll}
R_1 ^{p^3} = & R_1 \\
R_2 ^{p^3} = & Vect_{\mathbb{K}}\left\{  (x_i \cdot x_j) \cdot x_k-x_i \cdot (x_j \cdot x_k)+(x_j \cdot x_i) 
\cdot x_k-x_j \cdot (x_i \cdot x_k)  \right\}\\
R_3 ^{p^3} = & Vect_{\mathbb{K}}\left\{  (x_i \cdot x_j) \cdot x_k-x_i \cdot (x_j \cdot x_k)+(x_i \cdot x_k) 
\cdot x_j-x_i \cdot (x_k \cdot x_j)  \right\} \\
R_4 ^{p^3} = & Vect_{\mathbb{K}}\left\{  (x_i \cdot x_j) \cdot x_k-x_i \cdot (x_j \cdot x_k)+(x_k \cdot x_j) 
\cdot x_i-x_k \cdot (x_j \cdot x_i) \right\}\\
R_5 ^{p^3} = & R_5 \\
R_6 ^{p^3} = & Vect_{\mathbb{K}} \{  (x_1 \cdot x_2) \cdot x_3-x_1 \cdot (x_2 \cdot x_3) +
(x_2 \cdot x_3) \cdot x_1-x_2 \cdot (x_3 \cdot x_1) +(x_3 \cdot x_1) \cdot x_2\\                                          
&-x_3 \cdot (x_1 \cdot x_2) 
  +(x_2 \cdot x_1) \cdot x_3-x_2 \cdot (x_1 \cdot x_3) +(x_3 \cdot x_2) 
\cdot x_1-x_3 \cdot (x_2 \cdot x_1)\\
&   +(x_1 \cdot x_3) \cdot x_2-x_1 \cdot (x_3 \cdot x_2) \}
\end{array}$$
This corresponds to $R_i^{p ^3}=\K(\mathcal{O}(Y_i))$ with
$$Y_i=\Sigma _{\sigma  \in G_i} \sigma (X)$$
with $X=x_1\cdot (x_2\cdot x_3)-(x_1\cdot x_2)\cdot x_3.$ 

\noindent We denote by $(G_i\!-\!p^3\mathcal{A}ss)$ the corresponding quadratic operads. The corresponding dual 
operads are discribed by the following ideals of relations:
$$\begin{array}{ll}
(R_1 ^{p^3})^! = & R_1 \\
(R_2 ^{p^3})^! = & Vect_{\mathbb{K}}\left\{  (x_i \cdot x_j) \cdot x_k-x_i \cdot (x_j \cdot x_k);(x_i \cdot x_j) 
\cdot x_k+ (x_j \cdot x_i)\cdot  x_k \right\}\\
(R_3 ^{p^3})^! = & Vect_{\mathbb{K}}\left\{  (x_i \cdot x_j) \cdot x_k-x_i \cdot (x_j \cdot x_k);
(x_i \cdot x_j )\cdot x_k)+(x_i \cdot x_k) \cdot x_j  \right\} \\
(R_4 ^{p^3})^! = & Vect_{\mathbb{K}}\left\{  (x_i \cdot x_j) \cdot x_k-x_i \cdot (x_j \cdot x_k);
(x_i \cdot x_j) \cdot x_k+(x_k \cdot x_j) \cdot x_i \right\}\\
(R_5 ^{p^3})^! = & R_5^! \\
(R_6 ^{p^3})^! = & Vect_{\mathbb{K}} \{  (x_i \cdot x_j) \cdot x_k-x_i \cdot (x_j \cdot x_k) ;
(x_i \cdot x_j) \cdot x_k +(x_j \cdot x_i) \cdot x_k\\                                          
&(x_i \cdot x_j) \cdot x_k+(x_i \cdot x_k)\cdot x_j \}.
\end{array}$$
The proof is analogous to the Lie-admissible case. Let us note that these operads are also of rank $2$ except for $i=1$.

\medskip

\noindent {\bf 3. The  $\mathbb{K}\, [\Sigma_3]$-associative algebras}

\medskip

\noindent This example of nonassociative algebras generalizes the previous, considering not only the one dimensional
invariant subspace of $\mathbb{K}\, [\Sigma_3]$ but all the invariant subspaces.
Recall that, for every $v \in \mathbb{K}\,[\Sigma_3]$, we have denoted by $\mathcal{O}(v)$ 
  the corresponding orbit and by $F_v=\mathbb{K}\,(\mathcal{O}(v))$ the 
linear subspace
of $\mathbb{K}\,[\Sigma_3]$ generated by $\mathcal{O}(v).$ Since $F_v$ is a 
$\Sigma_3$-invariant subspace of $\mathbb{K}\,[\Sigma_3]$,  by Mashke's theorem,  
it is a direct sum of irreducible invariant subspaces.
Moreover,  given  an invariant subspace $F$ of $\mathbb{K}\,[\Sigma_3]$, there
exists $v \in \mathbb{K}\, [\Sigma_3]$ (not necessarily unique) such that $F=F_v=\mathbb{K}(\mathcal{O}(v)).$

\begin{definition}

1. A $\K$-algebra $(\mathcal{A},\mu)$ is called $\mathbb{K}\,[\Sigma_3]$-associative if there exists
$v \in \mathbb{K}\,[\Sigma_3], v \neq 0,$ such that 
$$A(\mu) \circ \Phi_v^\mathcal{A}=0$$
where $A(\mu )$ is the associator of $\mu $.

2. Let $ A(\mu) \circ \Phi_v^\mathcal{A} =0$ and $ A(\mu) \circ \Phi_w^\mathcal{A} =0$ be two identities 
satisfied by the algebra $(\mathcal{A},\mu)$. 
We say that these identities are equivalent if
$F_v= \mathbb{K}( \mathcal{O}(v))=F_w=\mathbb{K} (\mathcal{O}(w)).$ 
\end{definition}
Remark that if $F_v$ is not an irreducible 
invariant subspace, then there exists
$w \in F_v$ such that $F_w \subset F_v$ and $F_w \neq F_v.$ In this case the identity
$A(\mu) \circ \Phi_v ^\mathcal{A}=0$ implies $A(\mu) \circ \Phi_w^\mathcal{A} =0$ but these identities are not equivalent.

\medskip

\noindent {\bf Examples.}

\noindent 1. If $v=Id-\tau_{23}$, the relation 
\begin{eqnarray}
\label{2}
A(\mu) \circ \Phi_v^{\A}=0
\end{eqnarray}
becomes
$$ A(\mu)(x_1 \otimes x_2 \otimes x_3 -x_1 \otimes x_3 \otimes x_2)=0.
$$
The corresponding algebra is a pre-Lie algebra.

\noindent 2. The Lie-admissible and third-power associative algebras are $\mathbb{K}\,[\Sigma_3]$-associative algebras. 
In fact an algebra $(\A,\mu)$ is  Lie-admissible if
$A(\mu) \circ \Phi_V^\A =0$ and third-power associative if
$A(\mu) \circ \Phi_W^\A =0$ with 
$$V=Id-\tau _{12} -\tau _{23}-\tau _{13}+c_1+c_2$$
and
$$W=Id+\tau _{12} +\tau _{23}+\tau _{13}+c_1+c_2.$$

\noindent 3. For $i=1,...,6$ we denote by $V_i$ and $W_i$ the vectors of $\K[\Sigma _3]$ given by 
$$V_i=\sum_{\sigma \in G_i}(-1)^{\epsilon(\sigma)}\sigma, \quad W_i=\sum_{\sigma \in G_i}\sigma.$$
Then    $(\A,\mu)$ is a $G_i$-$\mathcal{A}ss$-algebra if
$$A(\mu) \circ \Phi_{V_i}^{\A}=0$$  
and a $G_i$-$p^3\mathcal{A}ss$-algebra if
$$A(\mu) \circ \Phi_{W_i}^\mathcal{A} =0.$$ 

\medskip

They are particular cases of $\K[\Sigma _3]$-associative algebras. 
We shall return, in the last section, on the determination of the corresponding operads.

\section{The operad $\tilde {\cal{P}}$ associated to a quadratic operad ${\cal{P}}$}

In the previous sections we saw that for some quadratic operads, the dual operad gives a way
to construct on the tensor product $ \A \otimes \B$ of a  $\cal{P}$-algebra $\A$ and a ${\cal{P}}^!$-algebra $\B$ 
a structure of $\cal{P}$-algebra for the usual tensor product $\mu_{\A \otimes \B}=\mu_{\A} \otimes \mu_{\B}$. 
But this is not true for every quadratic operad. 
In this section we define from a given quadratic operad ${\cal{P}}$ an associated quadratic operad, 
denoted by $\tilde {\cal{P}}$, whose fondamental property is to satisfy the above property on the  
tensor product. 

\medskip

Let $(\A,\mu)$ be a $\mathcal{P}$-algebra where $\mathcal{P}$ is a quadratic operad. Let $R$ be the submodule 
of $\Gamma (E)(3)$ defining the relations of $\A$. We denote by $A^L(\mu)=\mu \circ (\mu \otimes Id)$ and 
$A^R(\mu)=\mu \circ (Id \otimes \mu).$ If we suppose that $R$ is of rank $k$, 
then the multiplication $\mu$ satisfies $k$ relations of type
$$A^L(\mu)\circ \Phi_{v_i}^\A-A^R(\mu)\circ \Phi_{w_i}^\A=0$$
where  $v_i,w_i \in \mathbb{K}[\Sigma_3]$ for any $i \in I=\{1,..,k\}$ and the vectors $(v_i)_{i \in I}$ are linearly independent
as well as the vectors $(w_i)_{i \in I}$.

\medskip

\noindent {\bf Examples.} 

\medskip

\noindent 1. The associative algebras correspond to $k=1$ and $v_1=w_1=Id$, pre-Lie algebras to $k=1$
and $v_1=w_1=Id -\tau_{23}$ and more generally 
$G_i$-associative algebras correpond to $k=1$ and $v_1=w_1=V_i$.
The Lie-admissible algebras correpond to $k=1$, 
$v_1=w_1=V$ and the $3$-power associative algebras to $v_1=w_1=W$.

\medskip

\noindent 2. The Leibniz algebras correspond to $k=1$ and $v_1=Id-\tau_{23}$, $w_1=Id$.

\bigskip

\noindent Let $\mathcal{P}$ be a quadratic operad generated by $E \subset \mathbb{K}[\Sigma_2]$. For every 
$v=\Sigma_{l=1}^6S a_l\sigma _l \in \K[\Sigma _3]$ we consider on 
$\Gamma (E)(3)$ the linear maps
$$\Psi ^L_v ((x_i\cdot x_j)\cdot x_k) =\Sigma a_l\sigma _l((x_i\cdot x_j)\cdot x_k), \ 
\Psi ^L_v (x_i\cdot (x_j\cdot x_k))=0$$
and
$$\Psi ^R_v ((x_i\cdot x_j)\cdot x_k) =0, \ 
\Psi ^R_v (x_i\cdot (x_j\cdot x_k))=\Sigma a_l\sigma _l(x_i\cdot (x_j\cdot x_k)). $$
Let $R$ be the module of relations  of $\mathcal{P}$. If it is of rank $k$, it is written 
$$R= Vect_{\mathbb{K}}\left\{ 
 ( \Psi^L_{v_p}((x_i\cdot x_j)\cdot x_k)  - \Psi^R_{w_p}(x_1 \cdot ( x_2 \cdot x_3)), \ 1\leq p\leq k \right\} $$
with $v_p=\sum_{i=1}^6 a_i^p \sigma_i$ and 
$w_p=\sum_{i=1}^6 b_i^p \sigma_i$ for $1 \leq p \leq k.$ 

\noindent Let $\tilde{E}$ be the sub-module of $\mathbb{K}[\Sigma_2]$ defined by
$$\tilde{E}=\left\{
\begin{array}{l}
E \ \ \mbox{\rm{if} } E=1 \! \! 1 \oplus Sgn_2 \\
\mathcal{C}om(2)=1 \! \! 1 \ \ \mbox{\rm{if} } \ E=1 \! \! 1 \ \  \mbox{\rm{or} } \  Sgn_2 
\end{array}
\right.$$ 

\noindent If $dim \tilde E=2$, we denote
$\tilde{R}$ the $\mathbb{K}[\Sigma_3]$-module generated by the vectors 

$$
\left\{
\begin{array}{l}
a^{p}_i  a^{p}_j \Phi^L_{\sigma_i-\sigma_j} ((x_1\cdot x_2)\cdot  x_3), \\
\\
b^{p}_i  b^{p}_j \Phi^R_{\sigma_i-\sigma_j}(x_1 \cdot (x_2\cdot x_3)) , \\\\

a^{p}_i  b^{p}_j( \Phi^L_{\sigma_i}((x_1\cdot x_2)\cdot  x_3)   - \Phi^R_{\sigma_j}(x_1 \cdot( x_2\cdot x_3)), 
\end{array}
\right.$$
for $1\leq p\leq k$.

\noindent If $\tilde E=\mathcal{C}omm(2)$, $\tilde R$ is generated also by these vectors, but modulo the relations
of commutation. 

\begin{definition}
The operad $\tilde{\mathcal{P}}$ associated to the quadratic operad $\mathcal{P}$ is the quadratic operad 
generated by $\tilde{E}$ and with $\K[\Sigma_3]$-submodule of relations $\tilde{R}$.
\end{definition}
We have the main result :

\begin{theorem}
Let $\A$ be a $\mathcal{P}$-algebra and $\B$ a $\tilde{\mathcal{P}}$-algebra. Then  the algebra
$\A \otimes \B$ with product $\mu_{\A \otimes \B}$ is a $\mathcal{P}$-algebra.
\end{theorem}
{\it Proof. } Let $\A$ be a $\mathcal{P}$-algebra. Its multiplication $\mu_{\A}$ satisfies 
$$(A^L({\mu_{\A}}) \circ \Phi ^{\A}_{v_p}-A^R({\mu_{\A}})\circ \Phi ^{\A}_{w_p})(x_1\otimes x_2\otimes x_3)=0$$
for $p=1,...,k$. If $\B$ is a $\tilde{\mathcal{P}}$-algebra,  its multiplication $\mu_{\B}$ satisfies
$$
\left\{
\begin{array}{l}
A^L({\mu_{\B}}) \circ \Phi ^{\B}_{\sigma _i-\sigma _j}(y_1\otimes y_2\otimes  y_3) = 0, \ 
\mbox{\rm if} \ \exists p,  \ a^i_pa^j_p\neq 0 \\
\\
A^R({\mu_{\B}})\circ  \Phi ^{\B}_{\sigma _i-\sigma _j}(y_1\otimes  y_2\otimes  y_3) = 0, \ \mbox{\rm if} \ 
\exists p,  \ b^i_pb^j_p\neq 0 \\
\\
A^L({\mu_{\A}}) \circ \Phi ^{\B}_{\sigma _i}-A^R({\mu_{\B}})\circ  \Phi ^{\B}_{\sigma _j}(y_1\otimes  y_2\otimes  y_3)=0,
\ \mbox{\rm if} \  \exists p,  \ a^i_pb^j_p\neq 0. 
\end{array}
\right.
$$
Now we consider the product $\mu_{\A\otimes \B}$. We have
$$
\begin{array}{ll}
\medskip
&(A^L(\mu_{\A\otimes \B}) \circ \Phi ^{\AB}_{v_p}-A^R(\mu_{\A\otimes \B})\circ \Phi ^{\AB}_{w_p})
(x_1\otimes y_1\otimes x_2\otimes y_2\otimes x_3 \otimes y_3)\\
\medskip 
=&(\Sigma a_i^p A^L(\mu_{\A\otimes \B})\circ  \sigma _i
-\Sigma b_i^p A^R(\mu_{\A\otimes \B})\circ  \sigma _i)
(x_1\otimes y_1\otimes x_2\otimes y_2\otimes x_3 \otimes y_3) \\
\medskip 
=&\Sigma a_i^p(A^L(\mu_{\A})\circ  \sigma _i(x_1\otimes x_2\otimes x_3)\otimes 
A^L(\mu_{\B})\circ  \sigma _i(y_1 \otimes y_2 \otimes y_3)) \\
\medskip
& -\Sigma b_i^p(A^R(\mu_{\A})\circ  \sigma _i(x_1\otimes x_2\otimes x_3)\otimes 
A^R(\mu_{\B})\circ  \sigma _i(y_1 \otimes y_2 \otimes y_3)) \\
\medskip
=&(\Sigma a_i^p(A^L(\mu_{\A})\circ  \sigma _i(x_1\otimes x_2\otimes x_3))\otimes 
A^L(\mu_{\B})\circ  \sigma _j(y_1 \otimes y_2 \otimes y_3) \\
\medskip
& -(\Sigma b_i^p(A^R(\mu_{\A})\circ  \sigma _i(x_1\otimes x_2\otimes x_3))\otimes 
A^R(\mu_{\B})\circ  \sigma _j(y_1 \otimes y_2 \otimes y_3)) \\
\end{array}
$$
where $j$ is choosen in $\{1,\cdots 6\}$ such that $a_j^p \neq 0$,
$$
\begin{array}{ll}
=&\Sigma a_i^pA^L(\mu_{\A})\circ  \sigma _i(x_1\otimes x_2\otimes x_3) -
\Sigma b_i^pA^R(\mu_{\A})\circ  \sigma _i(x_1\otimes x_2\otimes x_3)) \\
\medskip
& \otimes 
A^R(\mu_{\B})\circ  \sigma _j(y_1 \otimes y_2 \otimes y_3)) \\
\medskip
=&(A^L(\mu_{\A}) \circ \Phi ^{\A}_{v_p}-A^R(\mu_{\A})\circ \Phi ^{\A}_{w_p})
(x_1\otimes x_2\otimes x_3)) 
 \otimes 
A^R(\mu_{\B})\circ  \sigma _j(y_1 \otimes y_2 \otimes y_3)) \\
\medskip
=& 0.
\end{array}
$$

\section{Some examples}

\subsection{$\mathcal{P}={G_i\!-\!\mathcal{A}ss}$}
\begin{proposition}
If $\mathcal{P}$ is a $(G_i\!-\!\mathcal{A}ss)$ operad, then the operads $\mathcal{P}^!$ and 
$\tilde{\mathcal{P}}$ are equal.
\end{proposition}
{\it Proof. }In this case we have $k=1$ and $v_1=w_1=V_i$.  Then $\tilde{\mathcal{P}}$ is defined by the module of relations
$$
\left\{
\begin{array}{l}
a_i  a_j \Psi^L_{\sigma_i-\sigma_j} ((x_1\cdot x_2)\cdot  x_3), \\
\\
a_i  a_j \Psi^R_{\sigma_i-\sigma_j}(x_1 \cdot (x_2\cdot x_3)) , \\\\

a_i  a_j( \Psi^L_{\sigma_i}((x_1\cdot x_2)\cdot  x_3)   - \Psi^R_{\sigma_j}(x_1 \cdot( x_2\cdot x_3)), 
\end{array}
\right.$$
where $V_i=\Sigma a_j\sigma _j$. This system is reduced to
$$
\left\{
\begin{array}{l}
 \Psi^L_{Id}((x_1\cdot x_2)\cdot  x_3)   - \Psi^R_{Id}(x_1 \cdot( x_2\cdot x_3)), 
a_i  a_j \Psi^L_{\sigma_i-\sigma_j} ((x_1\cdot x_2)\cdot  x_3), \\
\end{array}
\right.$$
which corresponds to the dual operad.

\noindent In particular, if $ \mathcal{P} =\mathcal{L}ie\mathcal{A}dm$, 
then $\mathcal{P}=\tilde{\mathcal{P}}= 
\mathcal{C}omm3$  where $\mathcal{C}omm3$ is the quadratic operad defined from the submodule of relations :
$$R=Vect_\mathbb{K}\{ (x_i \cdot x_j) \cdot x_k-x_i \cdot (x_j \cdot x_k), 
(x_i \cdot x_j) \cdot x_k-(x_j \cdot x_i) \cdot x_k \}.$$ 
Thus, a $\mathcal{C}omm3$-algebra $\A$ is $3$-commutative if it is associative and satisfies
$$x_i\cdot x_j \cdot x_k=x_{\sigma(i)} \cdot x_{\sigma(j)}  \cdot x_{\sigma(k)} $$
for every $\sigma \in \Sigma_3.$
If $\A$ is unitary this implies that $\A$ is a commutative algebra. 
If not, we have that $\A^2$ is contained in the center of $\A$. The associated Lie algebra is two step nilpotent.

\subsection{$\mathcal{P}=\mathcal{L}ie$ }
 If $\mathcal{P}=\mathcal{L}ie$ then $\tilde{\mathcal{P}}={\mathcal{P}}^!=\mathcal{C}om$. In fact, in this case,
$k=1$ and $v_1=w_1=Id+c_1+c_2$. As $v_1=w_1$, a $\tilde{\mathcal{P}}$-algebra is associative and the module of relation
corresponds to the operad $\mathcal{C}om$. 

\subsection{$\mathcal{P}=\mathcal{L}eib$}  
Let $\mathcal{P}=\mathcal{L}eib$ be the Leibniz operad. A Leibniz algebra is defined by the relation
$$x(yz)-(xy)z+(xz)y=0.$$
In this case the associated  $\tilde{\mathcal{P}}$ operad  corresponds to the relations
$$x(yz)=(xy)z$$
and
$$(xy)z=(xz)y.$$
Thus a $\tilde {\mathcal{L}eib}$-algebra is an associative algebra satisfying
$$xyz=xzy.$$
This relation is equivalent to
$$x[y,z]=0$$
with $[y,z]=yz-zy$. This last identity implies that if $x$ et $y$ are in the derived  Lie subalgebra of 
the associated Lie algebra, then $xy=0$. The derived Lie algebra is then abelian and the Lie algebra is
2 step nilpotent. The dual operad, also denoted by $\mathcal{Z}inb$, corresponds to the identity
$$(xy)z-x(yz)-x(zy)=0.$$ Thus a $\tilde {\mathcal{L}eib}$-algebra is a Zinbiel algebra (i.e. a 
$\mathcal{L}eib^!$-algebra) if $x(yz)=(xy)z=0$ (every product of $3$ elements of the associative algebra is zero). 
These algebras are  nilalgebras $\A $ satisfying $ \A^3=0$.
For example, any associative commutative  algebra is a $\tilde {\mathcal{L}eib}$-algebra. Every
$\tilde {\mathcal{L}eib}$-algebra with unit is commutative. In dimension $3$ the algebra defined by
$$e_1e_1=e_2, \ e_1e_3=e_3e_3=e_2$$
is a noncommutative $\tilde {\mathcal{L}eib}$-algebra.

\subsection{Determination of $\tilde{ \mathcal{P}}$ in the Lie-admissible case}

In the following table we describe the  multiplications of the algebras corresponding  
to the operads $\mathcal{P},\mathcal{P}^!,\tilde {\mathcal{P}}$ 
and this for the operads described in \cite{G.R2}, Theorem 3. 
Recall that these algebras are Lie-admissible  $\mathbb{K}\,[\Sigma_3]$-associatives algebras.
We denote by $A(x,y,z)$
the associator : $A(x,y,z)=(x\cdot y)\cdot z-x\cdot (y\cdot z)$.

\medskip

\noindent
$
\left\{
\begin{array}{l}
 \mathcal{P}=\mathcal{L}ie\mathcal{A}dm : A(x,y,z)-A(y,x,z)-A(x,z,y)-A(z,y,x)+ A(y,z,x)+A(z,x,y)=0.\\
 \mathcal{P}^!: A(x,y,z)=0, \ x \cdot y\cdot z=y\cdot x\cdot z=x\cdot z\cdot y. \\
\tilde {\mathcal{P}}=\mathcal{P}^! 
\end{array}
\right.
$

\medskip

\noindent 
$
\left\{
\begin{array}{l}
 \mathcal{P}=G_5-\mathcal{A}ss : A(x,y,z)+A(y,z,x)+A(z,x,y)=0.\\
 \mathcal{P}^!: A(x,y,z)=0, \ x\cdot y\cdot z=y\cdot z\cdot x=z\cdot x\cdot y. \\
\tilde {\mathcal{P}}=\mathcal{P}^! 
\end{array}
\right.
$

\medskip

\noindent 
$
\left\{
\begin{array}{l}
 \mathcal{P} : \alpha A(x,y,z)-\alpha A(y,x,z)+(\alpha+\beta-3)A(z,y,x)-\beta A(x,z,y)
 +\beta A(y,z,x)\\
\ \ +(3-\alpha-\beta)A(z,x,y)=0, \ (\alpha,\beta) \neq (1,1) \\
 \mathcal{P}^!: A(x,y,z)=0, (\alpha-\beta)(x\cdot y\cdot z-y\cdot x\cdot z)+
(\alpha+2\beta-3)(z\cdot y\cdot x-z\cdot x\cdot y)=0 \\
\end{array}
\right.
$

The computation of $\tilde {\mathcal{P}}$ depends on the values of the parameters $\alpha $ et $\beta $.  
If $(\alpha ,\beta ) \neq (3,0)$ ou $(0,3)$ or $(0,0)$ then 
$$\tilde {\mathcal{P}}=\mathcal{L}ie\mathcal{A}dm^! $$
If $(\alpha ,\beta )= (3,0)$ then $ \mathcal{P} = G_2-\mathcal{A}ss$ and 
$$\tilde {\mathcal{P}}= \mathcal{P}^!=G_2-\mathcal{A}ss^!.
$$
If $(\alpha ,\beta )= (0,3)$ then $ \mathcal{P} = G_4-\mathcal{A}ss $ and 
$$\tilde {\mathcal{P}}= \mathcal{P}^!=G_4-\mathcal{A}ss^!.
$$
If  $(\alpha ,\beta )=(0,0)$, then
$ \mathcal{P} = G_3-\mathcal{A}ss $ and
$$\tilde {\mathcal{P}}= \mathcal{P}^!=G_3-\mathcal{A}ss ^!.$$

\medskip

\noindent 
$
\left\{
\begin{array}{l}
 \mathcal{P} : A(x,y,z)+(1+t)A(y,x,z)+A(z,y,x)
+A(y,z,x)+(1-t)A(z,x,y)=0 , t \neq 1\\
 \mathcal{P}^!: A(x,y,z)=0,(t-1) x\cdot y\cdot z-(t-1) y\cdot x\cdot z-(t+2) z\cdot y\cdot x+(1+2t) x\cdot z\cdot y \\
\ \ -(1+2t)y\cdot z\cdot x+ (t+2) z\cdot x\cdot y=0.\\
\tilde {\mathcal{P}}=\mathcal{L}ie\mathcal{A}dm^!. 
\end{array}
\right.
$

\medskip

\noindent 
$
\left\{
\begin{array}{l}
 \mathcal{P} : 2A(x,y,z)+A(y,x,z)+A(x,z,y)
 +A(y,z,x)+A(z,x,y)=0. \\
 \mathcal{P}^!: A(x,y,z)=0,x\cdot y\cdot z+y\cdot x\cdot z-z\cdot y\cdot x-z\cdot x\cdot y=0 .\\
\tilde {\mathcal{P}}=\mathcal{L}ie\mathcal{A}dm^!. 
\end{array}
\right.
$

\medskip

\noindent 
$
\left\{
\begin{array}{l}
 \mathcal{P} : 2A(x,y,z)-A(y,x,z)-A(z,y,x)
 -A(x,z,y)+A(y,z,x)=0. \\
 \mathcal{P}^!: A(x,y,z)=0,x\cdot y\cdot z-y\cdot x\cdot z-z\cdot y\cdot x-x\cdot z\cdot y+y\cdot z\cdot x+z\cdot x\cdot y=0\\
\tilde {\mathcal{P}}=\mathcal{L}ie\mathcal{A}dm^!. 
\end{array}
\right.
$

\medskip

\noindent 
$
\left\{
\begin{array}{l}
 \mathcal{P}=\mathcal{A}ss : A(x,y,z)=0. \\
 \mathcal{P}^!= \mathcal{A}ss \\
\tilde {\mathcal{P}}=\mathcal{A}ss 
\end{array}
\right.
$
\begin{proposition}
Let $ \mathcal{P}$ be an operad corresponding to a  $\K [\Sigma _3]$-associative Lie-admissible algebra type.
Then $\tilde{ \mathcal{P}}= \mathcal{P}^!$ if and only if $ \mathcal{P}$ is the operad 
$G_i-{\mathcal{A}}ss$ for some  $i.$
\end{proposition}

We can also find the operads such that  $\tilde{ \mathcal{P}}= \mathcal{P}^!$ in the case of a quadratic 
operad generated by a commutative operation (i.e $E=1\! \! 1$) or an anticommutative one
(i.e $E=Sgn_2$)

\begin{proposition}
Let $ \mathcal{P}=\mathcal{P}(\K,E, R)$ be a quadratic operad generated by an operation with a symmetry
 (i.e. $E=Sgn_2 $ or $E=1\! \! 1$). Then $\tilde{ \mathcal{P}}= \mathcal{P}^!$ if and only if $E= Sgn_2 $  and $ \mathcal{P}=\mathcal{C}om$ or
$ \mathcal{P}=\Gamma(Sgn_2)$, the free operad generated by the signum representation.

\end{proposition}

\subsection{$\mathcal{P}=\mathcal{P}oiss$}

A Poisson algebra over $\mathbb{K}$ is a $\mathbb{K}$-vector space 
equipped with two bilinear products:

\noindent 1) a Lie algebra multiplication, denoted by $\{ , \}$, called the Poisson bracket,

\noindent 2) an associative commutative multiplication, denoted by
$\bullet$.

\smallskip

\noindent These two operations satisfy the Leibniz condition:
 \begin{eqnarray}
 \label{leib}
\{X   \bullet Y,Z \}=X\bullet \{Y,Z \}+\{ X,Z\}\bullet Y,
 \end{eqnarray}
for all $X,Y,Z.$
In \cite{MR}, one proves that a Poisson algebra can be defined by only
one nonassociative product, denoted by $X\cdot Y$, satisfying the following identity \begin{eqnarray}
\label{associator} 3A_{\cdot}(X,Y,Z)=(X\cdot Z)\cdot Y+(Y\cdot
Z)\cdot X-(Y\cdot X)\cdot Z-(Z\cdot X)\cdot Y,
\end{eqnarray}
where $A_{\cdot}(X,Y,Z)=(X\cdot Y)\cdot Z-X\cdot(Y\cdot Z)$ is the associator of the product $\cdot.$
The corresponding quadratic operad is with one generating operation and of rank $1$. Let us denote by 
$\mathcal{P}oiss$ this operad. If 
 $$ \Phi^L_{v_1}((x_i\cdot x_j)\cdot x_k)  - \Phi^R_{w_1}(x_1 \cdot ( x_2 \cdot x_3)) $$
is the generator of the module of relations $R$ of $\mathcal{P}oiss$, we have
$$v_1=3Id- \tau _{23}-c_1+\tau _{12}+c_2$$
and
$$w_1=3Id.$$
Then $\tilde{\mathcal{P}oiss}$ is generated by
$$
\left\{
\begin{array}{l}
(x_1\cdot x_2)\cdot x_3 - (x_1\cdot x_3)\cdot x_2\\
(x_1\cdot x_2)\cdot x_3  - (x_2\cdot x_3)\cdot x_1 \\
(x_1\cdot x_2)\cdot x_3  - (x_2\cdot x_1)\cdot x_3 \\
(x_1\cdot x_2)\cdot x_3 - (x_3\cdot x_1)\cdot x_2\\
(x_1\cdot x_2)\cdot x_3 - x_1\cdot (x_2\cdot x_3)
\end{array}
\right.
$$
and  $\tilde{ \mathcal{P}oiss}= \mathcal{C}omm3.$ 

\medskip

\noindent{\bf Remark. } In this work, we have considered only the classical product on the tensor product of algebras. It is possible 
to do the same work considering generalized or twisted tensor product. In these case we can probably 
define also an associated operad.

\noindent Let us consider two Poisson algebras $(\A, \mu _{\A} )$ and  $(\B,\mu _{\B})$ two
Poisson algebras defined by the nonassociative multiplication (\ref{associator}).  Let $\tau $ be the twist map:
$$
\tau(x\otimes y)=y\otimes x.
$$   
If we consider on $\A\otimes \B$ the following product
$$ \mu _{\A} \otimes _\tau \mu _{\B}=3\mu _{\A}\otimes \mu _{\B} -\mu _{\A}\otimes (\mu _{\B}\circ \tau )-
(\mu _{\A} \circ \tau )\otimes \mu _{\B}+ (\mu _{\A} \circ \tau )\otimes (\mu _{\B}\circ \tau )$$
then $(\A\otimes \B,\mu _{\A} \otimes _\tau \mu _{\B})$ is a Poisson algebra.

\end{document}